\documentclass[11pt]{article}

\input{amssym.def}
\input{amssym}
\usepackage{eufrak}

\newtheorem{theorem}{Theorem}

\newtheorem{definition}[theorem]{Definition}

\newtheorem{proposition}[theorem]{Proposition}
\newtheorem{remark}[theorem]{Remark}

\def\pa{{\partial}}
\def\u{\mbox{{\sf u}}}

\def\HH{{\cal H}}

\def\u{{\mathrm{u}}}
\def\v{{\mathrm{v}}}

\def\gg{{\frak{g}}}
\def\ov{{\overline}}
\def\Y{\mathbf{Y}}

\def\qq{q^{-1}}

\def\LL{{\cal{L}}}

\def\B{{\cal B}}

\def\de{\delta}

\def\ot{\otimes}

\def\C{{\Bbb C}}

\def\Z{{\Bbb Z}}
\def\Sym{{\rm Sym\, }}

\def\vv{V^{\otimes 2}}

\def\ov{\overline}

\def\al{{\alpha}}

\def\be{\begin{equation}}
\def\ee{\end{equation}}

\def\bea{\begin{eqnarray}}
\def\eea{\end{eqnarray}}

\begin{document}

\makeatletter
\renewcommand{\theequation}{{\thesection}.{\arabic{equation}}}
\@addtoreset{equation}{section} \makeatother

\title{Determinants in Quantum Matrix Algebras and Integrable Systems}

\author{\rule{0pt}{7mm} Dimitri
Gurevich\thanks{gurevich@ihes.fr}\\
{\small\it Universit\'e Polytechnique Hauts-de-France,  LMI}\\
{\small\it F-59313 Valenciennes, France}\\
{\small \it and}\\
{\small \it Interdisciplinary Scientific Center J.-V.Poncelet}\\
{\small\it Moscow 119002, Russian Federation}\\
\rule{0pt}{7mm} Pavel Saponov\thanks{Pavel.Saponov@ihep.ru}\\
{\small\it
National Research University Higher School of Economics,}\\
{\small\it 20 Myasnitskaya Ulitsa, Moscow 101000, Russian Federation}\\
{\small \it and}\\
{\small \it
Institute for High Energy Physics, NRC "Kurchatov Institute"}\\
{\small \it Protvino 142281, Russian Federation}}

\maketitle

\begin{abstract}
We define quantum determinants in  Quantum Matrix Algebras, related to couples of compatible braidings following the scheme from \cite{G}.  We
establish relations between these determinants and the so-called  column-(row-)determinants, often used in the theory of   integrable systems.
Also, we generalize the quantum integrable spin systems from \cite{CFRS} by using  generalized Yangians, related to couples of compatible braidings.
We demonstrate that  such quantum integrable spin systems are not uniquely  determined by  the "quantum coordinate ring" of the basic space $V$.  
For instance, the "quantum plane" $xy=qyx$ gives rise to two different integrable systems: rational and trigonometric ones.
\end{abstract}

{\bf AMS Mathematics Subject Classification, 2010:} 81R50, 81R12


\section{Introduction}

The notion of quantum (or $q$-)determinant was  introduced  in the papers of L.D. Faddeev's school (see, for instance, \cite{KS}) in connection with Quantum Inverse 
Scattering Method. Initially, such determinants were introduces in RTT algebras, associated with the $U_q(sl(N))$ $R$-matrices or with some current 
(i.e. depending on spectral parameters) $R$-matrices. 

Nevertheless,  in \cite{G} a large family of other involutive and Hecke symmetries was constructed and quantum determinants were defined in RTT algebras, associated 
with {\em even}\footnote{The term {\em even} means that the $R$-skew-symmetric algebra $\mathsf{\Lambda}_R(V)$ has a finite number of non-trivial homogenous 
components and  the highest non-trivial component $\mathsf{\Lambda}_R^m(V)$ is one-dimensional. In this case we say that $R$ is of rank $m$.} symmetries.

Let us recall that these symmetries are particular cases of braidings, that it operators $R:\vv\to\vv$ subject to the so-called {\em braid relation}
$$
(R\ot I)(I\ot R)(R\ot I)=(I\ot R)(R\ot I)(I\ot R).
$$
Hereafter, $V$ is a finite dimensional vector space (called basic) and $I$ stands for  the identity operator or matrix\footnote{Note that the operators $PR$ where $P$ is 
the flip, are subject to the so-called Quantum Yang-Baxter Equation and are usually called $R$-matrices.}.

A braiding $R$ is called a Hecke symmetry (resp.,  an involutive symmetry), if it meets a supplementary condition
$$ 
(R - qI)(R + \qq I)=0,\quad q\not=\pm 1,\qquad ({\rm resp.,}\,\,\, R^2=I).
$$
The nonzero parameter $q\in \C$ is assumed to be generic, i.e. such that
$$
k_q=\frac{q^k-q^{-k}}{q-\qq}\not=0,\quad \forall \,k\in \Z.
$$

The objective of the present paper is three-fold. First, by using the scheme from \cite{G}, we  define  quantum determinants  in  quantum matrix algebras (QMA), 
associated with couples of compatible braidings introduced in \cite{IOP1}, half-quantum algebras (HQA) defined in \cite{IO}, and generalized Yangians introduced 
in \cite{GS1, GS2}. Also, we show that in some RTT algebras  (including those associated with the quantum group (QG) $U_q(sl(N))$) quantum determinants  
can be cast under the form of column-(or row-)determinants, which are very popular in the literature on integrable systems.

Second, by using quantum elementary symmetric polynomials, closely related to  the quantum determinants, we exhibit the Bethe subalgebras in all generalized Yangians
in question. Thus, we get  quantum integrable systems, which are a far-reaching generalization  of the  spin systems from \cite{CFRS} and their rational counterparts. 
Different forms of the corresponding determinants are also discussed.

Third, we would like to draw the reader's attention to the fact that, contrary to a popular belief, the "quantum coordinate ring" of the basic space $V$ does not uniquely 
determine the corresponding quantum algebra and quantum determinant. Thus, the so-called "quantum plane" defined by $x y-q  y x=0$, gives rise to two completely 
different RTT algebras and generalized Yangians and consequently, leads  to two different integrable systems.

The paper is organized as follows. In the next section, by using the method of \cite{G}, we define the quantum determinants in the QMA associated with couples $(R,F)$ 
of compatible braidings, where $R$ is an even involutive or  Hecke symmetry. In section 3, we exhibit some properties of   the quantum determinants in different algebras.
In particular, we describe symmetries which allow defining column-(or row-)determinants. Also, we consider the quantum determinants in the right and left half-quantum
algebras. In Section 4 we define the quantum determinant in rational and trigonometric generalized Yangians and exhibit integrable systems associated 
with such Yangians.  

\medskip
\noindent
 {\bf Acknowledgement}  The authors are indebted to Vladimir Rubtsov for elucidating discussions. The work of P.S. was partially supported  by the RFBR grant 
19-01-00726-a.
 
\medskip
\noindent
{\bf  Conflict of Interest:} The authors declare that they have no conflicts of interest.

 \section{Quantum determinants in QMA}

With any Hecke symmetry\footnote{In this section we mainly deal with Hecke symmetries. The corresponding results and formulae for involutive symmetries can be 
obtained by setting $q=1$.} $R$ we associate the $R$-symmetric and $R$-skew-symmetric algebras of the space $V$ defined respectively as the following quotients
 of the free tensor algebra $T(V)$ of the space $V$ :
$$
{\rm Sym}_R(V)=T(V)/\langle {\rm Im}(qI-R)\rangle,\qquad \mathsf{\Lambda}_R(V)=T(V)/\langle {\rm Im}(\qq I+R)\rangle.
$$
Here $\langle J\rangle$ stands for a two-sided ideal generated by a set $J\subset T(V)$. The ground field is $\C$.

Each homogenous component ${\rm Sym}^{(k)}_R(V)$ (resp., $\mathsf{\Lambda}^{(k)}_R(V)$) can be identified with the image of the $R$-symmetrizer $S^{(k)}(R)$ 
(resp., $R$-skew-symmetrizer $A^{(k)}(R)$)  acting on the space $V^{\ot k}$. These projectors  can be defined by the following recursion formulae
\begin{eqnarray}
S^{(1)}=I,&\quad& S^{(k)}=\frac{1}{k_q} \,S^{(k-1)}\left(q^{-(k-1)} I+(k-1)_q R_{k-1\,k}\right)S^{(k-1)},\quad k\ge 2,\nonumber\\
 &&\label{q-symm}\\
A^{(1)}=I,&\quad& A^{(k)}=\frac{1}{k_q} \,A^{(k-1)}\left(q^{k-1} I-(k-1)_q R_{k-1\,k}\right)A^{(k-1)},\quad k\ge 2.\nonumber
\end{eqnarray}
As usual, the bottom indices indicate the positions  where  matrices or operators are located. These formulae  can be deduced from the representation theory of the 
symmetric group, provided  $R$ is involutive, or that of the Hecke algebra, provided $R$ is a Hecke symmetry (see \cite{Gy}).

Let $R$ and $F$ be braidings. Following \cite{IOP1} we say that the ordered couple $(R,F)$ is  compatible (or braidings $R$ and $F$ are compatible) if the following 
relations take place
$$
R_{12}\, F_{23}\, F_{12} = F_{23}\, F_{12}\, R_{23},\qquad R_{23}\, F_{12}\, F_{23} = F_{12}\, F_{23}\, R_{12}.
$$
Below, we always assume $R$ to be an involutive or a Hecke symmetry.

Let $L=\|l_j^i\|_{1\leq i,j\leq N}$ be an $N\times N$ matrix and $L_1=L\ot I_{2...p}$, $p\geq 2$. (Thus, $L_1$ is an $N^p\times N^p$ matrix.) We introduce the following 
notation:
\be
L_{\ov 1}=L_1,\quad L_{\ov{k+1}}=F_{k\,k+1} L_{\ov{k}}F_{k\,k+1}^{-1},\quad k\leq p-1.
\label{F-copy}
\ee
In the case $F=P$, where $P$ stands for the   usual flip, we recover the standard definition: $L_{k+1}=P_{k\,k+1}L_{{k}} P_{k\,k+1}$.

Now, we define a QMA $\LL(R,F)$ as a unital associative algebra generated by matrix elements of the matrix $L=\|l_i^j\|$ subject to the system of commutation relations:
\be
R_{12}L_{\ov 1} L_{\ov 2}=L_{\ov 1} L_{\ov 2} R_{12}.
\label{defin}
\ee
The matrix $L$ is called the {\em generating matrix} of the algebra $\LL(R,F)$.

Observe that the compatibility of the braidings $R$ and $F$ entails that the defining relations of the algebra $\LL(R,F)$ can be pushed forward to higher positions in the 
following sense
$$
R_{k\,k+1} L_{\ov k} L_{\ov{k+1}}=L_{\ov k} L_{\ov{k+1}} R_{k\,k+1},\quad \forall\,k< p.
$$

Note that each of the couples $(R,P)$ and $(R,R)$ is evidently compatible. The corresponding algebras $\LL(R,P)$ and $\LL(R,R)$ are respectively the RTT algebra and
Reflection Equation (RE) one\footnote{One more example of compatible braidings $(R,F)$ is that formed by braidings from (\ref{matrr}), where the second matrix plays the 
role of $R$. Also, an example of such couples is exhibited in \cite{IOP2}.}. The defining relations of the former algebra $\LL(R,P)$ read
$$
R_{12} L_{1} L_{2}=L_{ 1} L_{2} R_{12}.
$$
The defining relations of the RE algebra $\LL(R,R)$ can be cast under the following form:
$$
R_{12}L_{1}R_{12} L_{1}=L_{1}R_{12}L_{1}R_{12}.
$$

\begin{remark}\rm
It should be emphasized that if a symmetry $R$ is a deformation of the usual flip  the corresponding RTT and RE algebras are deformations of the commutative algebra 
$\Sym(gl(N))$, i.e. dimensions of the homogenous components of these QMA are classical (if $R$ is a Hecke symmetry, the parameter $q$ must be generic).

However, if $R$ is a braiding coming from the quantum groups of the series $B_n, C_n, D_n$, this property fails. Thus, any similar deformation of the algebra $\Sym(\gg)$, 
where $\gg$ is a Lie algebra belonging to one of these series, does not exist. By contrary, there exist quantum deformations of  the function algebra  $Fun(G)$, where
$G$ is the corresponding Lie group. The corresponding quotients of the RTT algebras are exhibited in \cite{FRT}.
\end{remark}

Now, assume the symmetry $R$ to be {\em even}. Let $\mathsf{\Lambda}_R^{(m)}(V),$ $m\geq 2$ be the highest non-trivial homogenous component\footnote{In general, 
$m$ could be different from $N=\dim V$ (see \cite{G,GPS2}).} of the algebra $\mathsf{\Lambda}_R(V)$. Since, by definition, the dimension  of this component is 1, then 
there exist two tensors
\be
\u=\|u_{i_1...i_m}\|\quad {\rm and}\quad \v=\|v^{j_1...j_m}\|, \label{tens} \ee
such that
$$
A^{(m)}(x_{i_1}\ot...\ot x_{i_m})=u_{i_1...i_m}\,v^{j_1...j_m}x_{j_1}\ot...\ot x_{j_m},\quad \langle \v,\u\rangle:=v^{i_1...i_m}\,u_{i_1...i_m}=1.
$$

Hereafter, $\{x_i\}_{1\leq i\leq N}$ is a basis of the space $V$ and summation over repeated indices is always understood. Thus, the element $v^{j_1...j_m}x_{j_1}\ot...\ot x_{j_m}$
is a generator of ${\rm Im}(A^{(m)})$. Note that the tensors $\u$ and $\v$ are defined up to a rescaling\footnote{If the rank of a symmetry $R$ equals 2, it is possible to recover 
the symmetry $R$ by knowing $\u$ and $\v$. All couples $(\u,\v)$, giving rise to such symmetries are classified in \cite{G}.}
$$
\u\mapsto a\u,\quad \v\mapsto a^{-1} \v,\qquad a\in \C,\,\,a\not=0.
$$

By analogy with \cite{G} we introduce the following definition.

\begin{definition} \rm
The element of QMA $\LL(R,F)$
\be
{\det}_{\LL(R,F)}(L):= \langle \v|L_{\ov 1}...L_{\ov m}|\u\rangle:= v^{i_1...i_m} \, (L_{\ov 1}...L_{\ov m})_{i_1...i_m}^{j_1...j_m} \,u_{j_1...j_m} ,
\label{det}
\ee
is called {\it the quantum determinant} of the generating matrix $L$.
\end{definition}

Of course, the quantum determinant  $\det_{\LL(R,F)}(L)$ can be written in other explicit forms, modulo the defining relations of the QMA $\LL(R,F)$. Some of them are 
presented below.

The form (\ref{det})  will be called {\em canonical}.

Consider now two examples. Let us fix a basis $\{x=x_1, y=x_2\}$ of the basic space $V$, $N= \dim \, V=2$ and introduce two symmetries, represented  by the following 
matrices in this basis
\be
\left(\begin{array}{cccc}
1&0&0&0\\
0&0&q&0\\
0&\qq&0&0\\
0&0&0&1
\end{array}\right),
\qquad
\left(\begin{array}{cccc}
q&0&0&0\\
0&q-\qq&1&0\\
0&1&0&0\\
0&0&0&q
\end{array}\right).
\label{matrr}
\ee
Each of these symmetries is a deformation of the usual flip $P$. The former symmetry is involutive, the latter one is a Hecke symmetry, coming from the QG $U_q(sl(2))$.

For the involutive symmetry we have
$$
\u=(u_{11},u_{12},u_{21},u_{22})=\frac{1}{2} (0, 1, -\qq, 0),\quad  \v=(v^{11},v^{12},v^{21},v^{22})=(0,1,-q,0).
$$
For the Hecke symmetry we have
$$
\u=\frac{\qq}{2_q} (0, 1, -q, 0),\quad  \v=(0,1,-q,0).
$$
Note that the tensors $\v$ corresponding to these symmetries coincide with each other and, consequently, the algebras 
\be
{\Sym}_R(V)=T(V)/\langle \v\rangle=T(V)/\langle xy-qyx \rangle,
\label{q-plane}
\ee
called the "quantum plane", are the same for the both symmetries  exhibited in (\ref{matrr}). Nevertheless, the tensors $\u$ are different. Consequently, the canonical forms 
of the corresponding determinants $\det_{\LL(R,F)}(L)$ differ from each other for all couples $(R,F)$.

Let us compute  these determinants for the corresponding RTT algebras $\LL(R,P)$ and  for the RE algebras $\LL(R,R)$. Denote entries of  the generating matrix $L$ of 
these algebras as follows $l_1^1=a$, $l_1^2=b$, $l_2^1=c$ and $l_2^2=d$:
$$
L=\left(\!
\begin{array}{cc}
a&b\\
c&d
\end{array}
\!\right).  \label{mattrr}
$$

\noindent
{\bf Example 1.}
The defining relations of the RTT algebra $\LL(R,P)$, corresponding to the first (involutive) matrix from (\ref{matrr})
$$
R_{12}L_1L_2 = L_1L_2R_{12}
$$
lead to the following systems for generators:
$$
ab=\qq ba,\quad ac=q ca,\quad ad=da,\quad bc=q^2 cb,\quad bd=q db,\quad cd=\qq dc.
$$

According to our definition, the canonical form of the quantum determinant in this algebra is
$$
{\det}_{\LL(R,P)}(L)=\frac{ad+da}{2} - \frac{\qq bc + qcb}{2}.
$$
With the use of the above commutation relations on the generators the canonical form can be transformed to the following expressions:
\be
{\det}_{\LL(R,P)}(L) = ad-q c b = a d-\qq bc.
\label{det1}
\ee

The defining relations between the generators of the algebra corresponding to the second matrix from (\ref{matrr}) are (see \cite{FRT})
\be
ab=q ba,\quad ac=q ca,\quad ad-da=(q-\qq)bc,\quad bc= cb,\quad bd=q db,\quad cd=q dc.
\label{RTT-rel}
\ee
The corresponding quantum determinant is
\be
{\det}_{\LL(R,P)}(L) =\frac{\qq a d+ q d a}{2_q}- \frac{bc + cb}{2_q} = ad-q cb
 =ad -q bc.
\label{det2}
\ee
The first expression in (\ref{det2}) is the canonical form of the determinant. Other expressions will be discussed in the next section. 

Let us also exhibit the corresponding algebras $\mathsf{\Lambda}_R(V)$. If $R$ is the first symmetry from (\ref{matrr}), then we have
$$
\mathsf{\Lambda}_R(V)=T(V)/\langle x^2, y^2, xy + q yx\rangle.
$$
If $R$ is the second symmetry from (\ref{matrr}), then the last generator of the ideal in the above quotient should be $qxy+yx$.

Thus, we see that the  quantum plane (\ref{q-plane}) gives rise to two different RTT algebras and consequently to two {\it different} determinants though we can find a specific 
form $ad-qcb$ which is the same for both determinants. In the next section we consider higher-dimensional analogs of these algebras and determinants in more details and 
explain this coincidence.

\medskip

\noindent
{\bf Example 2.} The defining relations
$$
R_{12}L_1R_{12}L_1 = L_1R_{12}L_1R_{12}
$$
of the RE algebra $\LL(R,R)$, corresponding to the first (involutive) matrix from (\ref{matrr}), in the explicit form read
$$
ab = ba, \quad ac = ca, \quad ad = da, \quad bc = cb, \quad bd = db, \quad cd = dc.
$$
Thus, this algebra is commutative.

It is not a surprising result since the involutive $R$ in (\ref{matrr}) is connected with the usual flip $P$ by the conjugation
$$
R_{12} = D_1P_{12}D_1^{-1}, \qquad
D=\left(\!
\begin{array}{cc}
q^{1/2}&0\\
0&q^{-1/2}
\end{array}
\!\right).
$$
Therefore, the matrix
$$
\tilde L = D^{-1}LD =
\left(\!
\begin{array}{cc}
a&b/q\\
qc&d
\end{array}
\!\right)
$$
satisfies the relation
$$
P_{12}\tilde L_1P_{12}\tilde L_1 = \tilde L_1P_{12}\tilde L_1 P_{12} \quad \Leftrightarrow\quad \tilde L_2\tilde L_1 = \tilde L_1\tilde L_2,
$$
which means that the entries of $\tilde L$ generate a commutative algebra. Therefore, the entries of the matrix $L$ also commute with each other.

The canonical form of the determinant in this case reads
$$
{\det}_{\LL(R,R)}(L) = \frac{ad+da}{2} - \frac{bc+cb}{2},
$$
or, taking into account the commutativity of the RE $\LL(R,R)$, it can be reduced to the classical expression $ad-bc$.

At last, if we take the second (Hecke) symmetry $R$ from (\ref{matrr}), we get the following system on the generators
\be
\begin{array}{lclcl}
q^2ab = ba&\quad &q^2ca = ac &\quad&ad = da\\
\rule{0pt}{5mm}
q(bc-cb) = \lambda a(d-a)&&q(cd - dc) = \lambda ca&&q(db-bd) = \lambda ab,
\end{array}
\label{REA-rel}
\ee
where we denoted $\lambda = q-q^{-1}$. The canonical form of the determinant is as follows
$$
{\det}_{\LL(R,R)}(L) = \frac{q(ad+da)}{2_q} - \frac{q(bc+q^2cb)}{2_q} - \frac{\lambda a^2}{2_q}.
$$
Upon taking into account the relations between the generators, we can transform the canonical determinant to equivalent forms:
\be
{\det}_{\LL(R,R)}(L) = ad - q^2cb = q^2(ad - bc) -q\lambda a^2.
\label{REA-det}
\ee

Another way of introducing quantum analogs of the determinant  is based on the notion of the quantum elementary symmetric polynomials defined via {\it the quantum traces}.
Such a quantum trace is well known in the cases related to the QG $U_q(sl(N))$.  Nevertheless, the quantum trace can be associated with any {\em skew-invertible} braiding $R$
by means of the following method belonging to V.Lyubashenko \cite{L1, L2}.

We say that a given braiding  $R:\vv\to\vv$ is {\em skew-invertible} if there exists an operator  $\Psi:\vv\to\vv$ such that
$$
{\rm Tr}_{(2)}R_{12}\Psi_{23}=P_{13} = {\rm Tr}_{(2)}\Psi_{12}R_{23} \quad \Leftrightarrow \quad R_{ij}^{kl}\Psi_{l m}^{jn}=\de_m^k \, \de_i^n =
\Psi_{ij}^{kl}R_{l m}^{jn}.
$$

If $R$ is a skew-invertible braiding, the corresponding {\it $R$-trace} ${\rm Tr}_R$ is defined by the formula
$$
{\rm Tr}_R X={\rm Tr} (C^RX),\quad C^R:={\rm Tr}_{(2)} \Psi.
$$
Here $X$ is an arbitrary $N\times N$ matrix (may be with noncommutative entries).

Consider a compatible couple of braidings $(R,F)$ and suppose the braiding $R$ to be skew-invertible. Now, we define a quantum version of the elementary symmetric 
polynomials in the algebra $\LL(R,F)$ as follows
\be
e_0=1,\quad e_k = {\rm Tr}_{R(1\dots k)}A^{(k)}L_{\overline 1}L_{\overline 2}\dots L_{\overline k},\quad k\ge 1.
\label{elem}
\ee
Hereafter, ${\rm Tr}_{(1\dots k)}={\rm Tr}_{(1)}\dots {\rm Tr}_{(k)}$.

By using the equality
\be
A^{(k)} L_{\ov 1}\dots L_{\ov{k}}=A^{(k)} L_{\ov 1}\dots L_{\ov{k}} A^{(k)}, 
\label{equa} 
\ee
valid in any QMA, we get the following relation
\be
{\rm Tr}_{R(1\dots m)} A^{(m)} L_{\ov 1}\dots L_{\ov{m}} ={\rm Tr}_{R(1\dots m)} A^{(m)} L_{\ov 1}\dots L_{\ov{m}} A^{(m)}=(\v\cdot_R \u)
\langle \v|L_{\ov 1}\dots L_{\ov m}|\u\rangle,
\label{dett}
\ee
where
\be
(\v\cdot_R \u)=v^{j_1...j_m}(C^{R})_{j_1}^{i_1}(C^{R})_{j_2}^{i_2}...(C^{R})_{j_m}^{i_m}\,   u_{i_1...i_m}.
\label{nota}
\ee
Thus, the highest elementary symmetric polynomial $e_m$ differs from the quantum determinant $\det_{\LL(R,F)}(L)$  by a numerical factor. In the particular case $F=P$ these
elements are just equal to each other since in this case $(\v\cdot_P \u)=1$ (note that $C^P=I$).

\section{Some properties of quantum determinants}

In this section we consider two questions. The first one is: what is the relation between the determinant ${\det}_{\LL(R,F)}(L)$ and the characteristic polynomial of the matrix $L$? 
The second question is: whether the quantum determinant is central? We always assume the rank of a symmetry $R$ to be $m$.

We say that a monic polynomial $ch(t)$ of degree $m$ is {\em characteristic} if $ch(L)=0$. In virtue of the Cayley-Hamilton theorem in the classical case $R=F=P$ (the 
corresponding algebra $\LL(P,P)$ is commutative) the characteristic polynomial reads
$$
ch(t) = {\det}_{\LL(P,P)}\,(L-t\, I).
$$

\begin{proposition} In the algebras $\LL(R,R)$, where $R$ is a Hecke symmetry, the following relation holds\rm
\be
{\det}_{\LL(R,R)}\,(L-t\, I)=\sum_{0\leq k\leq m} (-t)^{m-k} \al_k\, e_k,
\label{cha}
\ee\it
where $\displaystyle\al_k= q^{mk}\,\frac{m!}{k!(m-k)!}\,\frac{k_q!(m-k)_q!}{m_q!}$.
\end{proposition}

\noindent
{\bf Proof.} For any even Hecke symmetry of the rank $m$ the quantity (\ref{nota}) reads
$$
(v\cdot_R u)=q^{-m^2}.
$$
This follows from the relation (see \cite{GS2})
\be
{\rm Tr}_{R(k+1...m)}\, A^{(m)}_{1\dots m}= q^{-m(m-k)} \, \frac{k_q!(m-k)_q!}{m_q!} \,A^{(k)}_{1\dots k}
\label{chech}
\ee
for $k=0$.

Now, in the expansion of the element
$$
q^{m^2}{\det}_{\LL(R,R)}\,(L-t\, I)= {\rm Tr}_{R(1...m)}A^{(m)}_{1\dots m}(L-tI)_{\ov 1}...(L-t I)_{\ov m}
$$
we put together the terms containing $k$ factors $L_{\ov i}$ on some places and  the identity matrices on other positions. The number of such terms is
$\frac{m!}{k!(m-k)!}$ and they are equal to each other. This property is due to the fact that\footnote{Emphasize, that if $F\not=R$ this property fails.}
$$
{\rm Tr}_{R(1\dots m)}A^{(m)}_{1\dots m}\,L_{\ov i_1}L_{\ov i_2}\dots L_{\ov i_k}={\rm Tr}_{R(1\dots m)} A^{(m)}_{1\dots m}\,L_{\ov 1}L_{\ov 2}\dots L_{\ov k},
$$
for any ordered subset of indices $1\le i_1<i_2<\dots < i_k\le m.$ Now, it suffices to apply formula (\ref{chech}). Finally, we find
$$
{\rm Tr}_{R(1\dots m)}A^{(m)}_{1\dots m}\,L_{\ov 1}\dots L_{\ov k} = q^{-m(m-k)} \, \frac{k_q!(m-k)_q!}{m_q!}\,e_k(L).
$$
The proof is completed. \hfill \rule{6.5pt}{6.5pt}

\medskip

If $R$ is an involutive symmetry, then by setting $q=1$ in (\ref{cha}) we get the following claim.

\begin{proposition}
If $R$ is an involutive symmetry, the polynomial ${\det}_{\LL(R,R)}\,(L-t\, I)$ is characteristic.
\end{proposition}

If $R$ is a Hecke symmetry, the polynomial ${\det}_{\LL(R,R)}\,(L-t\, I)$ is {\it not}  characteristic. However, we get the characteristic polynomial upon replacing
$\al_k$ by $q^k$ in the right hand side of (\ref{cha}), that is we have:
$$
ch(t):=t^m-q  e_1\, t^{m-1}+q^2 e_2\, t^{m-2}+...+(-q)^{m-1}  e_{m-1}\, t+(-q)^{m} e_m.
$$
Thus, upon substituting $t=L$ in this polynomial, we get the Cayley-Hamilton identity for the matrix $L$:
$$
L^m-q  e_1\, L^{m-1} +q^2e_2\,  L^{m-2}+...+(-q)^{m-1}  e_{m-1}\, L+(-q)^{m} e_m\, I=0.
$$

Note that the first proof of this identity in the algebras $\LL(R,R)$ was given in \cite{GPS1}.

Now, we pass to the second question.  It is well known that if a Hecke symmetry $R$ comes from the quantum group $U_q(sl(N))$, the quantum determinant
${\det}_{\LL(R,P)}(L)$ is central (see \cite{FRT}). If in a given RTT algebra $\LL(R,P)$ the quantum determinant is central, then, by imposing the condition 
${\det}_{\LL(R,P)}(L)=1$, we can define a Hopf algebra structure in the quotient algebra.

However, in general the quantum determinant is {\it not} central in the algebras $\LL(R,P)$. As was shown in \cite{G}, the quantum determinant $\det_{\LL(R,P)}(L)$ is 
central if and only if  the matrix
$$
M=\|M_i^j\|,\quad \mathrm{where}\quad  M_i^j=u_{i\, i_2...i_m}\, v^{i_2... i_{m}\,j},
$$
is scalar.

Let us study  the centrality  of the quantum determinants in the algebras $\LL(R,P)$, corresponding to the symmetries (\ref{matrr}). By straightforward computations we get
respectively the matrices $M$ for the symmetries from (\ref{matrr}):
$$ 
- \frac{1}{2}
\left(\begin{array}{cc}
q&0\\
0&\qq
\end{array}\right)
\quad \mbox{and} \quad - \frac{1}{2_q}
\left(\begin{array}{cc}
1&0\\
0&1
\end{array}\right).
$$
Thus, the quantum determinant is not central in the algebra $\LL(R,P)$, corresponding to the involutive symmetry from (\ref{matrr})  and it is central in the algebra, corresponding 
to the Hecke symmetry from (\ref{matrr}).

By contrast, in the RE algebras $\LL(R,R)$ the quantum determinant is {\it always} central. Therefore, by imposing the condition  ${\det}_{\LL(R,R)}(L)=1$, we get a braided Hopf
algebra structure in the quotient algebra (see \cite{GPS2}).

Now,  we are going to discuss a way of reducing the quantum determinants to the so-called column-determinants and row-determinants,  playing an important role  in the theory of
integrable systems.

We set $F=P$, i.e. consider the RTT algebra $\LL(R,P)$. By using the relation (\ref{equa}) we get the following equality
$$
u_{i_1...i_m} v^{j_1...j_m} l_{j_1}^{k_1}...l_{j_m}^{k_m}= u_{i_1...i_m}\langle \v|L_1...L_m|\u\rangle v^{k_1...k_m}.
$$
Since the tensor $\u\not\equiv 0$, the factors $u_{i_1...i_m}$ can be cancelled. Assume also, that $m=N$ and $v^{12...N}=1$. This condition can be satisfied by a proper 
normalization of $\v$, provided  $v^{12...N}\not=0$. Then we get
\be 
{\det}_{\LL(R,P)}(L)=v^{j_1...j_N} l_{j_1}^{1}...l_{j_N}^{N}. 
\label{col} 
\ee

It is this form of the quantum determinant which is called {\it the column-determinant} of $L$. It is characterized by the property, that in each of its summands the multipliers 
$l_i^j$ are set  in the order of columns of the matrix $L$ enumerated  by the upper indices of $l_i^j$. Note that if $m\not= N$, we have no privileged component of the tensor $\v$ (like the 
component $v^{12...N}$).  

In a similar way, if the multipliers in each of the summands of the determinant are set in the order of rows of the matrix $L$, we call it {\it the row-determinant}. 
If $m=N$ and $u_{12...N}\not=0$, we can transform the canonical determinant $\det_{\LL(R,P)}$ to the form of row-determinant, which looks as follows
\be
{\det}_{\LL(R,P)}(L)= l_{1}^{i_1}...l_{N}^{i_N} u_{i_1...i_N}. 
\label{row} 
\ee

Observe that the column-determinant (resp., row-determinant) depends only on the tensor $\v$ (resp.,  $\u$). Thus, if two symmetries have the same tensors $\v$ but different
tensors $\u$ we have the identical form of the column-determinants but different row-determinants. It is just the case of the symmetries (\ref{matrr}) considered above. We see 
that the column-determinants (the middle expressions in (\ref{det1}) and (\ref{det2})) are equal to each other but the row-determinants (the right expressions) are diffident.

Now, introduce the higher-dimensional counterparts of the symmetries (\ref{matrr}) and exhibit the corresponding quantum determinants in the algebras $\LL(R,P)$.

The Hecke symmetry $R$ coming from the QG $U_q(sl(N))$ is as follows
$$
R_{ij}^{kl}=q^{\de_{k,l}}\de_j^k\,\de_i^l+(q-\qq) \theta_{(l>k)} \de_i^k\, \de_j^l,
$$
where $\theta_{(l>k)}=1$ if $l>k$  and $\theta_{(l>k)}=0$ if $l\le k$.  

Introduce, also, the involutive symmetry $R$ by its action on the basis vectors $x_i\otimes x_j$ of the space $V^{\otimes 2}$:
$$
R(x_i\ot x_i)=x_i\ot x_i,\quad R(x_i\ot x_j)=q\, x_j\ot x_i\,\,\,\mbox{\rm if}\,\, i<j,\quad R(x_i\ot x_j)=\qq\, x_j\ot x_i \,\,\,\mbox{\rm if}\,\, i>j.
$$

For the both symmetries the components of the tensors $\u$ and $\v$  are  nontrivial iff their indices are pairwise distinct. 
Then for the both symmetries the nontrivial components of the tensor $\v$ can be taken as
$$
v^{j_1... j_N}=(-q)^{l(\sigma)},
$$
where $l(\sigma)$ is the length (i.e. the minimal number of transpositions) of the permutation
$$
\sigma: (1...N)\mapsto (j_1...j_N).
$$
Sometimes, such a tensor $\v$ is called  the $q$-Levi-Civita tensor.

As for the tensors $\u$ for these symmetries, they are respectively equal to
$$
u_{i_1... i_N}=\alpha_1\, (-q)^{l(\sigma)},  \qquad u_{i_1... i_N}=\alpha_2\, (-\qq)^{l(\sigma)},  
$$
where 
$$
\al_1^{-1} = q^{\frac{N(N-1)}{2}}N_q!,\qquad \al_2^{-1} = N!
$$ 
are  normalizing factors.

Similarly to the 2-dimensional example above we have for both symmetries the same "quantum coordinate ring" 
$$ 
x_i\, x_j=q \,x_j\, x_i,\quad \forall\, i<j, 
$$
often called the {\it quantum torus} (with the additional condition that the generators are invertible). 
 
The formulae for the quantum column-determinants are also the same in both RTT algebras:
 $$
{ \det}_{\LL(R,P)}(L)=\sum_{\sigma}(-q)^{l(\sigma)} l_{j_1}^{1}...l_{j_N}^{N}.
$$
 In this form the quantum determinant $\det_{\LL(R,P)}(L)$ is given in \cite{FRT} for the RTT algebra associated with $U_q(sl(N))$ symmetry $R$.
 
 However, the  tensors $\u$ corresponding to the symmetries under consideration are different.  Consequently,  the row-determinants in the corresponding RTT algebras
 differ from each other. 
  
In conclusion of the section, we turn to the so-called half-quantum algebras (HQA) and the corresponding quantum determinants. 

Again, consider a compatible couple $(R,F)$, where $R$ is a skew-invertible Hecke symmetry. Introduce two systems of relations on the generating matrix $L=\|l^i_j\|_{1\leq i,j \leq N}$:
\be 
S^{(2)} L_{\ov 1} L_{\ov 2} A^{(2)} =0\quad \Leftrightarrow \quad L_{\ov 1} L_{\ov 2} A^{(2)}= A^{(2)} L_{\ov 1} L_{\ov 2} A^{(2)} , 
\label{first} 
\ee
\be 
A^{(2)} L_{\ov 1} L_{\ov 2}S^{(2)}=0\quad \Leftrightarrow \quad A^{(2)} L_{\ov 1} L_{\ov 2} = A^{(2)} L_{\ov 1} L_{\ov 2} A^{(2)},
\label{second}
\ee
where $R$-symmetrizer $S^{(2)}$ and $R$-skew-symmetrizer $A^{(2)}$ are defined in (\ref{q-symm}). The matrices with overlined  indices have the same meaning as above
(see (\ref{F-copy})).

The following claim is well known and can be checked straightforwardly.

\begin{proposition}
The system {\rm (\ref{defin})} is equivalent to the union of the systems {\rm (\ref{first})} and {\rm (\ref{second})}.
\end{proposition}

By imposing  on generators only {\it half} of the relations (only (\ref{first}) or (\ref{second})), we get a bigger algebra than that $\LL(R,F)$. Nevertheless, even in such an algebra
it is possible to develop some elements of linear algebra. We refer the reader to the paper \cite{IO}, where these algebras were introduced and studied. Some particular cases 
of these algebras were also considered   in \cite{CFR} and \cite{CFRS} under the name of {\em Manin matrices} and $q$-Manin-matrices.  

We call the right (resp., left) half-quantum algebra  (HQA) a unital algebra defined by the system (\ref{first})  (resp., (\ref{second})). They will be respectively denoted
 $\HH_r(R,F)$ and $\HH_l(R,F)$. 

If $R$ is an even symmetry, we  define the quantum determinant in the algebra $\HH_{\epsilon}(R,F)$, $\epsilon\in\{r,l\}$ by formula (\ref{det}) and denote it 
$\det_{\HH_{\epsilon}(R,F)}(L)$.

Quantum elementary symmetric polynomials can be still defined in the algebras $\HH_{\epsilon}(R,F)$ by formulae (\ref{elem}), where the projectors $A^{(k)}$ 
can be transferred to the rightmost position or put on the both positions: to the right and to the left of the chain of $L$-matrices..  

Note that the quantum determinant in the algebra $\HH_{\epsilon}(R,F)$ also differs from the highest quantum elementary symmetric polynomial by a factor.

\begin{remark} \rm Initially such type algebras were considered by Yu.Manin in \cite{M}. Their definition is motivated by the following consideration. Let us endow 
the space $V$ with the coaction of an RTT algebra: $x_i\to t_i^j\ot x_j$ and extend it to the space  $\vv$ by assuming that the generators $x_i$ and $t_k^j$  commute 
with each other. Then the relation (\ref{first}) (resp., (\ref{second})), where we put $F=P$, means that the subspace  ${\rm Im} \,A^{(2)}$  (resp., ${\rm Im} \,S^{(2)}$)  
is preserved under this coaction. However, if $F\not=P$, the assumption of commutativity of the generators $x_i$ and  $t_k^j$  is not suitable any more.
\end{remark}

Since the relation (\ref{equa}) is valid in any left HQA, the quantum determinant $\det_{\HH(R,P)_l}(L)$,  where $R$ is one of the symmetries (\ref{matrr}) or its 
higher-dimensional counterpart, can be cast under the form of column-determinant.  Whereas in any right HQA the quantum determinants can be cast under the form of
row-determinant. 

Nevertheless, in any HQA the number of relations between the generators is insufficient to prove that the elementary polynomials commute with each other.

\section{Generalized Yangians and  integrable systems of CFRS type}

 First, describe the Baxterization procedure which enables us to construct  current braidings via involutive and Hecke symmetries (see \cite{J,Jo}).
Let us precise that by a current  braiding $R(u,v)$ we mean  an operator depending on parameters and subject to the braid relation of the following form:
\be
R_{12}(u,v)R_{23}(u,w)R_{12}(v,w)=R_{23}(v,w)R_{12}(u,w)R_{23}(u,v).
\label{brel}
\ee

Given an involutive symmetry $R$, we associate with it a current braiding by the rule
\be
R(u,v)=R-\frac{I}{u-v},
\label{curr-inv}
\ee
whereas for a Hecke symmetry $R$ the corresponding current braiding reads
\be
R(u,v)=R-\frac{(q-\qq)u\, I}{u-v}.
\label{curr-Hecke}
\ee
By a straightforward calculation one can verify that these operators do satisfy the relation (\ref{brel}). The current braidings (\ref{curr-inv}) and (\ref{curr-Hecke}) (and all 
corresponding algebras) will be respectively called the {\em rational} and {\em trigonometric} ones.

Introduce a countable set of elements $l^i_j[k]$, $k\in\Z_{\ge 0}$, $1\le i,j\le N$, and consider a formal power series
\be
L(u)=\sum_{k\geq 0} L[k] u^{-k}, \quad L[k]=\|l_i^j[k]\|_{1\leq i,j\leq N}, 
\label{ser}
\ee
that is $L(u)$ is an $N\times N$ matrix and its entries are power  series in  $u^{-1}$ with coefficients $l^i_j[k]$.

A {\it generalized Yangian} $\Y(R,F)$ is an associative unital algebra generated by elements $l^i_j[k]$ subject to the system
\be
R_{12}(u,v)L_{\ov 1}(u) L_{\ov 2}(v)-L_{\ov 1}(v)L_{\ov 2}(u) R_{12}(u,v)=0,
\label{dl}
\ee
where $L_{\ov 1}(u)=L_{1}(u)$ and $L_{\ov 2}(u)=F_{12} L_{\ov 1}(u)F^{-1}_{12}$. Note that expanding the current matrix $L(u)$ is a series as indicated in  (\ref{ser}),
we get a countable set of {\it polynomial} relations on the generators $l_i^j[k]$.

In the case $F=R$ the algebra $\Y(R,R)$ with a supplementary condition $L[0]=I$ is called a {\em braided or generalized Yangian of RE type} (see \cite{GS1} for detail). Note 
that the condition $L[0]=I$ is motivated by the evaluation morphism, similar to that  in the Drinfeld's Yangian $\Y(gl(N))$. Note that the Drinfeld's Yangian  is a particular case 
of the rational Yangians, corresponding to the symmetry  $R=F=P$.

Observe that for a special value of the ratio $u/v=q^2$ in the trigonometric case the system (\ref{dl}) can be treated in terms of the HQA. A similar treatment is possible in
the rational case if $u-v=1$.

More precisely,  for the indicated relations between the parameters $u$ and $v$, the current braiding $R(u,v)$ becomes equal (up to a numerical factor) to the
$R$-skew-symmetrizer ${A}^{(2)}$. Thus, in the rational (respectively, trigonometric) case we arrive at the relations
\be 
A^{(2)}L_{\ov 1}(u)L_{\ov 2}(u-1)=L_{\ov 1}(u-1)L_{\ov 2}(u)A^{(2)}, 
\label{od} 
\ee
\be 
A^{(2)} L_{\ov 1}(u) L_{\ov 2}(q^{-2}\,u)= L_{\ov 1}(q^{-2}\,u) L_{\ov 2}(u) A^{(2)}.
\label{dv} 
\ee

Consequently, we have
$$
A^{(2)}L_{\ov 1}(u) L_{\ov 2}(u-1) S^{(2)}=0,\qquad  S^{(2)}L_{\ov 1}(u-1) L_{\ov 2}(u) A^{(2)}=0,
$$
in the rational case or, respectively,
$$
A^{(2)} L_{\ov 1}(u) L_{\ov 2}(q^{-2}\,u)S^{(2)}=0,\qquad  S^{(2)} L_{\ov 1}(q^{-2}\,u) L_{\ov 2}(u)A^{(2)}=0
$$
in the trigonometric case.

Using the formal Taylor series expansions
$$
L(u-1)= e^{-\pa_u}L(u)e^{\pa_u},\qquad L(q^{-2}u)= q^{-2u\pa_u}L(u)q^{2u\pa_u},
$$
where $\pa_u=\frac{d}{d\, u}$, we can rewrite the above relations in the form
$$
A^{(2)}(e^{-\pa_u} L_{\ov 1}(u)) (e^{-\pa_u}L_{\ov 2}(u))S^{(2)}=0,\qquad \ S^{(2)}(e^{\pa_u} L_{\ov 1}(u)) (e^{\pa_u}L_{\ov 2}(u))A^{(2)}=0
$$
in the rational case and 
$$
A^{(2)}(q^{-2u\pa_u} L_{\ov 1}(u))(q^{-2u\pa_u} L_{\ov 2}(u))S^{(2)}=0,\qquad  S^{(2)}(q^{2u\pa_u} L_{\ov 1}(u))(q^{2u\pa_u} L_{\ov 2}(u))A^{(2)}=0
$$
in the trigonometric case. 

Thus, the operator $e^{-\pa_u}L(u)$ (respectively, $q^{-2u\pa_u}\,L(u)$) plays the role of the generating matrix of a left HQA, whereas the operator $e^{\pa_u}L(u)$ 
(respectively, $q^{2u\pa_u}\,L(u)$) plays the role of the generating matrix of a right  HQA.

Now, define the quantum elementary symmetric polynomials in respectively rational and trigonometric cases as follows
$$
e_0(u)=1,\quad e_k(u) = {\rm Tr}_{R(1\dots k)}A^{(k)}\,L_{\overline 1}(u)\,L_{\overline 2}(u-1)\dots L_{\overline k}(u-k+1),\quad k\ge 1,
$$
$$
e_0(u)=1,\quad  e_k(u) = {\rm Tr}_{R(1\dots k)}A^{(k)}\,L_{\overline 1}(u)\,L_{\overline 2}(q^{-2}u)\dots L_{\overline k}(q^{-2(k-1)}u),\quad k\ge 1.
$$
Observe that in these formulae the projectors $A^{(k)}$ can be put after the chain of the matrices $L$ or at two positions: to the right and to the left of the chain of the matrices
$L$. Besides, the parameters of the matrices can be written in the inverse order. All these transformations lead to identical results.

Let $R$ to be an even symmetry of rank $m$ and define the quantum determinants in the generalized Yangians (respectively, rational and trigonometric) by setting
$$
{\det}_{\Y(R,F)} (L(u))=\langle \v|L_{\overline 1}(u)\,L_{\overline 2}(u-1)\dots L_{\overline m}(u-m+1)|\u\rangle,
$$
$$
{\det}_{\Y(R,F)} (L(u))=\langle \v|L_{\overline 1}(u)\,L_{\overline 2}(q^{-2}u)\dots L_{\overline m}(q^{-2(m-1)}u)|\u\rangle.
$$
Note, that the quantum determinant equals to the highest elementary symmetric polynomial $e_m(u)$ up to a factor in full analogy with (\ref{dett}). If $F=P$ this equality is exact. 

Consider the case $F=P$ in more detail. Assuming that $m=N$ and that $v^{1...N}$ and $u_{1...N}$ are non-trivial, we can cast the quantum determinant under the form of 
the column-determinant or the row-determinant:
$$ 
{\det}_{\Y(R,P)}(L)=v^{j_1...j_N} l_{j_1}^{1}(u)...l_{j_N}^{N}(u-N+1) = u_{i_1...i_N} l^{i_1}_{1}(u-N+1)...l^{i_N}_{N}(u)
$$
in the rational case and 
$$
 {\det}_{\Y(R,P)}(L)=v^{j_1...j_N} l_{j_1}^{1}(u)...l_{j_N}^{N}(q^{-2(N-1)}u) = u_{i_1...i_N} l^{i_1}_{1}(q^{-2(N-1)}u)...l^{i_N}_{N}(u)
$$
in the trigonometric one. 

Now,  let us take as $R$ the first (involutive) symmetry from (\ref{matrr}) or its higher-dimensional counterpart. Since the corresponding generalized Yangian is rational, we can 
write the corresponding quantum determinant as follows
$$
{\det}_{\Y(R,P)}(L)=\sum_{\sigma}(-q)^{l(\sigma)} l_{\sigma(1)}^{1}(u)...l_{\sigma(N)}^{N}(u-N+1) = \sum_{\sigma}(-q^{-1})^{l(\sigma)} l^{\sigma(1)}_{1}(u-N+1)...l^{\sigma(N)}_{N}(u).
$$

In the case of the $U_q(sl(N))$ symmetries $R$, the corresponding  generalized Yangian is trigonometric. Consequently, we have
$$
{\det}_{\Y(R,P)}(L)=\sum_{\sigma}(-q)^{l(\sigma)} l_{\sigma(1)}^{1}(u)...l_{\sigma(N)}^{N}(q^{-2(N-1)}u) = \sum_{\sigma}(-q)^{l(\sigma)}
l^{\sigma(1)}_{1}(q^{-2(N-1)}u)...l^{\sigma(N)}_{N}(u).
$$
We point out that the order of arguments in the matrices are opposite in the expressions for the column-determinant and for the row-determinant.

Also, note that similarly to the case of the algebras $\LL(R,F)$, the quantum determinant is always central in the generalized Yangians of  RE type $\Y(R,R)$ but it is not so
in these of RTT type $\Y(R,P)$. More precisely, the quantum determinant ${\det}_{\Y(R,P)}(L)$ is central iff it is so for the quantum determinant ${\det}_{\LL(R,P)}(L)$. This
property is proved in \cite{GS1}. 

The quantum versions of power sums can be also defined in all generalized Yangians  and some quantum versions of the Cayley-Hamilton identity are valid in all of them.

The subalgebra generated in the generalized Yangian $\Y(R,F)$ by the quantum elementary polynomials is called the Bethe subalgebra and denoted $\B(R,F)$.

\begin{proposition}{\rm  \cite{GSS}} 
For any compatible couple $(R,F)$ of braidings such that $R$ is a skew-invertible involutive or Hecke symmetry the Bethe subalgebra $\B(R,F)\subset \Y(R,F)$ is commutative.
\end{proposition}

A particular case of this claim, corresponding to $F=P$ and $R$ coming from the QG $U_q(sl(N))$ is proved in \cite{CFRS}. (Note that the formula for the projectors
$A^{(k)}$ should be taken as in (\ref{q-symm}))

The generalized Yangians of the RE type have a very important property: they admit evaluation morphisms. These morphisms were constructed in \cite{GS1}. Similarly to the 
evaluation map in the Drinfeld's Yangian $\Y(P,P)$ they have the form
$$
L(u)\mapsto 1+\frac{M}{u}
$$
both in the rational and trigonometric cases. However, the target algebras, generating by matrix elements of $M$, are different: for the rational generalized Yangians it is a
modified RE algebra, respective to the symmetriy $R$, and in the trigonometric case, it is a non-modified RE algebra.

In conclusion, we make a short remark. As we have already noticed, it is not possible to prove the commutativity of the quantum elementary polynomials in a HQA. However,
the relations (\ref{od}) and (\ref{dv}) enable one to establish this property in the generalized Yangians because they are more restrictive than the defining relations in the HQA.

\end{document}